\documentclass[12pt]{article}
\newlength{\abstractwidth}
\newcommand{\be}{\begin{equation}}
\newcommand{\bea}{\begin{eqnarray}}
\newcommand{\eea}{\end{eqnarray}}
\newcommand{\ee}{\end{equation}}

\newtheorem{theorem}{Theorem}

\begin{document}
\centerline{\bf \large THE ALMOST HERMITIAN-EINSTEIN FLOW}

\bigskip
\centerline{Chiung-Nan Tsai}
\bigskip

\centerline{Department of Mathematics}
\centerline{Columbia University, New York, NY 10027}
\bigskip
\bigskip

\begin{abstract}
{\small Let $X$ be a compact K\"ahler manifold, $E\rightarrow X$ a Hermitian
vector bundle and $L\rightarrow X$ an ample line bundle. We construct a non-linear heat flow
corresponding to the almost Hermitian-Einstein equation introduced by N.C. Leung,
and prove that the solution exists for a short time. We also construct a potential function
$D_k$ for this flow. In particular, $D_k$ decreases along the flow.}

\end{abstract}

\baselineskip=15pt

\section{Introduction and Main Results}

The problem of existence and uniqueness of Hermitian-Einstein metrics for Mumford-Takemoto stable
bundles has been solved by Uhlenbeck and Yau in \cite{UY} using the continuity method and
by Donaldson  in \cite{D} using the heat equation method. Under a smoothness
assumption on the bundles, Leung in \cite {L1} and \cite{L2} proved a similar result in the case
of almost Hermitian-Einstein metrics for Gieseker stable bundles, using a singular perturbation
technique and the result obtained by Donaldson-Uhlenbeck-Yau. This motivates us to study the almost
Hermitian-Einstein metrics by using the heat equation method.

\medskip

The setting is the following. Let $E\to X$ be a holomorphic vector bundle over an
$n$-dimensional K\"ahler manifold
$(X,\omega)$, and $L$ an ample line bundle over $X$
with $c_1(L)=[\omega]$. The bundle $E$ is called Hermitian
when it is equipped with a Hermitian metric $H=H_{\bar\alpha\beta}$.
Following Leung \cite{L1,L2}, the metric $H_{\bar\alpha\beta}$
is said to be
almost-Hermitian-Einstein if it satisfies for large positive $k$ the following
equation, called the almost-Hermitian-Einstein equation,
\be
\label{AHE}
[\exp({\frac{i}{2 \pi}F + k \omega I})Td(X)]^{(2n)} = \frac{1}{rk(E)}
\chi(X,E\otimes L^k)\frac{\omega ^n}{n!} I
\ee
where $I$ is the identity,
$F$ is the curvature 2-form on the Hermitian vector bundle $(E,H)$,
$Td(X)$ is the harmonic representative with respect to $\omega$ of the Todd class of $X$ ,
$rk(E)$ is the rank of $E$ and $\chi(X,E\otimes L^k)$ is the Euler characteristic of the bundle
$E\otimes L^k$:
$$
\chi(X,E\otimes L^k)=\sum^n_{i=0}(-1)^i {\rm dim}\, H^i(X,E\otimes L^k).$$
The leading terms in the almost Hermitian-Einstein equation
correspond to the coefficients of $k^n$, and keeping these terms only
gives the well-known Hermitian-Einstein equation
\be
\label{HE}
\Lambda F=\mu_E I,
\ee
where $\mu_E= c_1(E)[\omega]^{n-1}/{\rm rk}\,E$ is the slope of $E$,
and
$\Lambda F\equiv \hat F =g^{j\bar k}F_{\bar kj}$,
with the K\"ahler form given by $\omega={i\over 2}g_{\bar kj}dz^j\wedge d\bar z^k$.
The full almost Hermitian-Einstein equation
can be expressed in the following form
\be
\label{AHE2}
\frac{i}{2\pi} F \wedge \frac{{\omega}^{n-1}}{(n-1)!} =
 \mu_E \frac{{\omega}^n}{n!} + \sum_{j=1}^{n-1} \frac{1}{k^j} T_{j+1}
\ee
with the terms $T_j$ defined by
\be
T_j = \chi_E ^j \frac{{\omega}^n}{n!} - \sum_{k=0}^{j}
(\frac{i}{2\pi} F)^k Td_X ^{j-k}
\frac{{\omega}^{n-j}}{(n-j)!}.
\ee
For our purposes, it is convenient to rewrite (\ref{AHE}) also as
\be
\label{AHE1}
\Lambda F = \mu_{E} I - S(k),
\ee
where  $S(k)$ consists of terms involving powers
of $1/k$ strictly less than $n$, and is defined so that
the equation (\ref{AHE2}) coincide with the equation (\ref{AHE1}).

\medskip
Associated to the almost Hermitian-Einstein equation is the following natural flow
of endomorphisms $h(t)$ of $E$,
\be
\label{AHEheat}
\dot{h}(t)h^{-1}(t) = -(\ \Lambda F - \mu_E I + S(k)\ ),
\qquad h(0)=I,
\ee
where $H_0$ is a fixed Hermitian metric on $E$,
$H(t)$ is a Hermitian metric evolving with time $t$, and we have set
$h(t)=H(t) H_0^{-1}$.
The primary purpose of this paper is to study this flow, and in
particular, to construct a potential function for it.

\medskip

To construct a potential function, we consider {\it any} path $H(t)$ of Hermitian metrics
with $H(0)=H_0$, and
we introduce the following functional $D_k$ defined by
\be
\label{D}
D_k(H(t),H_0) = \int_{X}n R_2 \wedge {\omega}^{n-1} - \mu_E R_1 {\omega}^n + \int_0^t tr
(\sum_{j=1}^{n-1} {\frac{1}{k^j}} T_{j+1}\dot{h}h^{-1}).
\ee
where  $R_1=\log {\rm det}\,(tr(H H_0^{-1}))$, $R_2=\sqrt{-1}\int_0^t tr(F\dot h h^{-1})dt$ are
the well-known secondary characteristic classes.
Then our main results are as follows:

\bigskip

\begin{theorem}

{\rm (a)} The functional $D_k$ is a potential function
on the space of Hermitian metrics in the following sense.
If $h(s,t)$ satisfying
$h(s,0)=I$,
$h(s,1)=HH_0^{-1}$,
is a smooth deformation between the two paths
$h(0,t)H_0$ and $h(1,t)H_0$ joining two fixed
Hermitian metrics $H_0$ and $H$,
then $D_k(h(s,1)H_0,H_0)$ is independent of $s$.
Thus $D_k(h(s,1)H_0,H_0)$ can be considered as a function
of just the two end-points $H_0$ and $H$.

{\rm (b)} The functional $D_k$ is a Hamiltonian for the
flow (\ref{AHEheat}) in the sense that
\be
{d\over dt}D_k(H(t),H_0)
=
\int_X tr\,(\ (\Lambda F-\mu I+S(k))\ \dot h h^{-1})\,\omega^n.
\ee
In particular, $D_k$ is a decreasing function of $t$ along
the flow (\ref{AHEheat}).

{\rm (c)} For $k$ large, the flow
(\ref{AHEheat}) is parabolic, and hence exists for some
interval $0\leq t<T$.

\end{theorem}

\bigskip

\begin{theorem}

Let $\mu$ be the moment map, that is, the expression
defined by
the left hand side of the equation (\ref{AHE}). Then $\mu$ satisfies the
following evolution equation,
\be
\label{muflow}
\dot \mu =\sum_m \sum_{l=0}^{m-1}[\frac{i}{2\pi}C_l^{m-1}\omega ^l (\frac{i}
{2\pi}F)^{m-1-l} Td_{n-m}(X) \bar{\partial}\partial_H({n!\over\omega^n}\mu)]_{(sym)}k^l.
\ee
\end{theorem}

\bigskip
\noindent
{\it Acknowledgements.} I wolud like to thank Professor Mu-Tao Wang for many
discussions on this work. Special thanks are due to my advisor, Professor D.H. Phong,
for suggesting this problem to me and  for his constant encouragement. He helps me correct several mistakes and make this paper more readable.   This paper will be a part of the my future Ph.D. thesis in the Mathematics Department of Columbia University.

\section{Leung's work}

Almost Hermitian-Einstein metrics are natural from several viewpoints.
As shown by Leung,
their existence is equivalent to the notion of Gieseker stability.
They also have a natural interpretation in terms of symplectic geometry.
We give a brief review of some of Leung's
work in this section.

\medskip

Let the tangent vectors to the space of connections at $D_A$
be identified with the space of $End(E)$-valued one forms on $X$.
Then in \cite{L1,L2}, a one-parameter family of gauge-invariant 2-forms on the
space of
connections on E was defined as follows:
\be
\Omega_k (D_A)(B,C)=\int_X Tr_E [B\wedge \exp{(k\omega I+ \frac{i}{2\pi}F_A)}\wedge C]_{sym}Td(X)
\ee
where B and C are $End(E)$-valued one-forms on $X$,
and $F_A$ is the curvature tensor of the
connection $D_A$. Then, it can be shown that
$\Omega_k$ is a symplectic form, and the corresponding moment map
for the action of the gauge group is
\be
\label{moment}
\mu(D_A)=[{\rm exp}{(k\omega I +\frac{i}{2\pi}F)}Td_X]^{(2n)}.
\ee
Symplectic quotients can then be constructed as level sets of the moment map.
The almost Hermitian-Einstein equation
(\ref{AHE}) can be interpreted then as the equation for $\mu^{-1}(c)$,
where $c$ represents the right hand side of (\ref{AHE}).

\medskip

On the other hand, the notion of Gieseker stability of a vector bundle is the following.
Let E be a rank r holomorphic vector bundle (or coherent torsion-free sheaf in general)
over a projective variety X with ample line bundle L. The bundle $E$ is said to be
Gieseker stable if for any
nontrivial coherent subsheaf S of E, we have
\be
\frac{\chi(X,S\otimes L^k)}{rk\, S}< \frac{\chi(X,E\otimes L^k)}{rk\, E}
\ee
for large enough $k$. Then it is also shown in \cite{L1,L2}
that if $E$ is an irreducible sufficiently smooth holomorphic vector bundle over a compact
K\"ahler manifold X, then E is Gieseker stable if and only if there exists an almost
Hermitian-Einstein metric on E.

\section{Proof of Theorem 1}

We begin with the proof of (a). By a basic lemma of Donaldson,
the expression
$$
\int_{X}n R_2 \wedge {\omega}^{n-1} - \mu_E R_1 {\omega}^n
$$
is independent
of the paths of Hermitian metrics.
See \cite{S}, whose methods we adapt below.
Thus the two first expressions in the functional $D_k$ are path-independent,
and it suffices to show
that each single term in the remaining third expression
in (\ref{D}) is
also path-independent.
This means that we can simply deal with the following term:
\be
R_3 = \int_0^t \int_X tr(F^k Td^{j-k} \dot{h}h^{-1}{\omega}^{n-j}).
\ee
We begin by evaluating the second derivative $R_{3 t s}$ of $R_3$ with respect to
$r$ and $s$,
\bea
R_{3 t s} &=& \int_X tr \frac{d}{ds}(F^k Td^{j-k} h_t
h^{-1}{\omega}^{n-j}) \nonumber\\
&=&\int  tr (F_s F^{k-1} Td^{j-k} h_t h^{-1}{\omega}^{n-j})
+ \int tr (F F_s F^{k-2} Td^{j-k} h_t
h^{-1}{\omega}^{n-j})\nonumber\\
&&+ \int tr (F^2 F_s F^{k-3} Td^{j-k} h_th^{-1}{\omega}^{n-j})
+ \cdots \nonumber\\
&&+ \int tr (F^{k-1} F_s Td^{ j-k} h_t h^{-1}{\omega}^{n-j})
+ \int tr (F^{k} Td^{j-k} h_{ts} h^{-1}{\omega}^{n-j})\nonumber\\
&&- \int tr (F^k Td^{j-k} h_t h^{-1} h_s
h^{-1}{\omega}^{n-j}).
\eea
By the commutation formula $$F_s = \bar{\partial}\partial_H(h_s
h^{-1}) = - \partial_H\bar{\partial}(h_s h^{-1})- F h_s h^{-1} + h_s
h^{-1} F ,$$
we get
\begin{eqnarray}
R_{3 t s} &=& \int_X tr ((- \partial_H\bar{\partial}(h_s h^{-1})- F h_s
h^{-1} + \underline{h_s
h^{-1} F}) F^{k-1} Td^{j-k} h_t h^{-1}{\omega}^{n-j}) \nonumber\\
&&+ \int tr (F (- \partial_H\bar{\partial}(h_s h^{-1})- F h_s h^{-1} +
h_s
h^{-1} F) F^{k-2} Td^{j-k} h_t h^{-1}{\omega}^{n-j}) \nonumber\\
&&+ \int tr (F^2 (- \partial_H\bar{\partial}(h_s h^{-1})- F h_s h^{-1}
+ h_s
h^{-1} F) F^{k-3} Td^{j-k} h_th^{-1}{\omega}^{n-j})\nonumber\\
&&+ \cdots \nonumber\\
&&+ \int tr (F^{k-1} (- \partial_H\bar{\partial}(h_s h^{-1})- F h_s
h^{-1} + h_s
h^{-1} F) Td^{ j-k} h_t h^{-1}{\omega}^{n-j}) \nonumber\\
&&+ \int tr (F^{k} Td^{j-k} h_{ts} h^{-1}{\omega}^{n-j})
- \underline{\int_M tr (F^k Td^{j-k} h_t h^{-1} h_s
h^{-1}{\omega}^{n-j})}.\nonumber
\end{eqnarray}
Here, we can see the second term in each line and  the third term in the next line cancel.
The underlined terms also cancel, because $$tr AB = tr BA$$ if $A$ and $B$ are even-form-valued
matrices of the same size. Hence $R_{3ts}$ becomes
\bea
R_{3 t s} &=& \int_X tr ((- \partial_H\bar{\partial}(h_s h^{-1}))F^{k-1} Td^{j-k}
h_t h^{-1}{\omega}^{n-j}) \nonumber\\
&&+ \int_X tr (F (- \partial_H\bar{\partial}(h_s h^{-1}))F^{k-2} Td^{j-k} h_t h^{-1}{\omega}^{n-j})
\nonumber\\
&&+ \int_X tr (F^2 (- \partial_H\bar{\partial}(h_s h^{-1})) F^{k-3} Td^{j-k} h_th^{-1}{\omega}^{n-j})
\nonumber\\
&&+ \cdots \nonumber\\
&&+ \int_X tr (F^{k-1}  (-\partial_H\bar{\partial}(h_s h^{-1})) Td^{ j-k} h_t h^{-1}{\omega}^{n-j}) \nonumber\\
&&- \int_X tr (F^{k}  h_s h^{-1}  Td^{ j-k} h_t h^{-1}{\omega}^{n-j})
+ \int_X tr (F^{k} Td^{j-k} h_{ts} h^{-1}{\omega}^{n-j}).
\eea
Next, we compute $\int_X tr
\frac{d}{dt}(F^k Td^{j-k} h_s h^{-1}{\omega}^{n-j})$. By using this time
\be
F_t = \bar{\partial}\partial_H(h_t h^{-1})
\ee
we get
\bea
&&
\int_X tr
\frac{d}{dt}(F^k Td^{j-k} h_s h^{-1}{\omega}^{n-j})
\nonumber\\
&=&
\quad
\int_X tr(\bar{\partial}\partial_H(h_t h^{-1})F^{k-1} Td^{j-k} h_s
h^{-1}{\omega}^{n-j})
+ \int_X tr (F  \bar{\partial}\partial_H(h_t h^{-1})F^{k-2} Td^{j-k} h_s h^{-1}{\omega}^{n-j})
\nonumber\\
&&\quad
+ \int_X tr (F^2  \bar{\partial}\partial_H(h_t h^{-1}) F^{k-3} Td^{j-k} h_s h^{-1}{\omega}^{n-j})
+ \cdots \nonumber\\
&&
\quad+ \int_X tr (F^{k-1}  \bar{\partial}\partial_H(h_t h^{-1} ) Td^{ j-k} h_s h^{-1}{\omega}^{n-j})
+ \int_X tr (F^{k} Td^{j-k} h_{ts} h^{-1}{\omega}^{n-j})\nonumber\\
&&\quad
- \int_X tr (F^k Td^{j-k} h_s h^{-1} h_t h^{-1}{\omega}^{n-j}).
\eea
We can evaluate now  $ \int_X tr \frac{d}{ds}(F^k Td^{j-k} h_t
h^{-1}{\omega}^{n-j}) - \int_M tr \frac{d}{dt}(F^k Td^{j-k} h_s
h^{-1}{\omega}^{n-j}).$
We notice that the final two terms of $ \int_X tr \frac{d}{ds}(F^k Td^{j-k} h_t
h^{-1}{\omega}^{n-j})$ and $\int_M tr \frac{d}{dt}(F^k Td^{j-k} h_s
h^{-1}{\omega}^{n-j})$  cancel. Next, by using integration by parts, the Bianchi identity,
and the fact that Todd classes are
combinations of Chern classes, which are all closed forms, on the rest of the terms of
$ \int_X tr \frac{d}{ds}(F^k Td^{j-k} h_t
h^{-1}{\omega}^{n-j})-\int_M tr \frac{d}{dt}(F^k Td^{j-k} h_s
h^{-1}{\omega}^{n-j})$, we get
\begin{eqnarray}
&&\int_X tr \frac{d}{ds}(F^k Td^{j-k} h_t
h^{-1}{\omega}^{n-j}) - \int_X tr \frac{d}{dt}(F^k Td^{j-k} h_s
h^{-1}{\omega}^{n-j})\nonumber\\
&=& \int_X tr (( \bar{\partial}(h_s h^{-1})F^{k-1} Td^{j-k} \partial_H(h_t h^{-1}){\omega}^{n-j})) \nonumber\\
&&+ \int_X tr (F (\bar{\partial}(h_s h^{-1})F^{k-2} Td^{j-k} \partial_H(h_t h^{-1}){\omega}^{n-j}) )\nonumber\\
&&+ \cdots \nonumber\\
&&+ \int_X tr (F^{k-1} (\bar{\partial}(h_s h^{-1}) Td^{j-k} \partial_H(h_t h^{-1}){\omega}^{n-j}) )\nonumber\\
&&+ \int_X tr (( \partial_H(h_t h^{-1})F^{k-1} Td^{j-k} \bar{\partial}(h_s h^{-1}){\omega}^{n-j})) \nonumber\\
&&+ \int_X tr (F (\partial_H(h_s h^{-1}) F^{k-2} Td^{j-k} \bar{\partial} (h_sh^{-1}){\omega}^{n-j}))\nonumber\\
&&+ \cdots \nonumber\\
&&+ \int_X tr (F^{k-1} (\partial_H(h_s h^{-1})Td^{j-k} \bar{\partial} (h_sh^{-1}){\omega}^{n-j}))\nonumber
\end{eqnarray}
We claim the sum of these $2k$ terms is zero. We consider the $i$-th term and the
$(2k-i+1)$-th term as a pair. For simplicity, assume $i=1$, then the pair
\be
\int_X tr ( \bar{\partial}(h_s h^{-1})F^{k-1} Td^{j-k}
\partial_H(h_t h^{-1}){\omega}^{n-j})
\ {\rm and}
\ \int_X tr (F^{k-1} (\partial_H(h_s h^{-1})Td^{j-k} \bar{\partial} (h_sh^{-1}){\omega}^{n-j}))
\nonumber
\ee
cancel, because $tr AB= -tr BA$ if $A$ and $B$ are odd-form-valued matrices of the same size.
The same argument works on the other pairs, and this proves the claim.
Let $h(s,t)$ be now a deformation of paths as in the statement
of Theorem 1. Then we have
\bea
\frac{d}{ds} R_3 &=& \frac{d}{ds}\int_{t=0}^1 \int _X tr(F^k h_t
h^{-1} Td^{j-k}{\omega}^{n-j}) \nonumber \\
&=& \int_{t=0} ^1 \int_X \frac{d}{dt}tr(F^k h_s
h^{-1}Td^{j-k}{\omega}^{n-j}) \nonumber\\
&=& \int_{t=0} ^1 \int_X \frac{d}{dt}tr(F^k h_s
h^{-1}Td^{j-k}{\omega}^{n-j}) \nonumber\\
&=& \int_X tr(F^k h_s h^{-1}Td^{j-k}{\omega}^{n-j})|_{t=0} ^{t=1}
=0,
\eea
since $h_s=0$ at $t=0$ and $1$. This finishes the proof of (a) in  Theorem 1.

\medskip
Next, we prove (b). By the work of Donaldson, we already know that
${\partial\over\partial t}R_1=tr(\dot h h^{-1})$ and
${\partial\over\partial t}R_2=\sqrt{-1}tr(F\dot h h^{-1})$.
It follows immediately that
\bea
\frac{dD_k}{dt} &=& \int_X tr((\Lambda F - \mu_E I) \dot{h} h^{-1}) {\omega}^{n} +
tr(\sum_{j=1}^{n-1} \frac{1}{k^j} T_{j+1}\dot{h}h^{-1}) \nonumber\\
              &=& \int_X tr(\ (\Lambda F - \mu_E I + S(k)\ )\dot{h}
h^{-1})\,{\omega}^{n}.
\eea
In particular, along the flow (\ref{AHEheat}), we have
\bea
\frac{dD_k}{dt} = -\int_X |(\Lambda F - \mu_E I + S(k))|^2{\omega}^{n}
\leq 0,
\eea
and (b) is proved.

\medskip

Finally, (c) follows from the fact
that the principal symbol of $-\wedge F$ is actually the
same as the symbol of the Laplacian $\Delta$. Right away we can see the
principal symbol of the left hand side of almost Hermitian-Einstein flow is
elliptic at the initial time. The short-time existence is then a consequence
of the general theory of parabolic equations. Q.E.D.

\section{Proof of Theorem 2}

To prove Theorem 2, we write the almost Hermitian-Einstein flow in
the following form:
\bea
\dot{h}h^{-1} = \frac{-[\exp({i\over 2 \pi}F + k \omega I)Td(X)]^{(2n)} +
\frac{1}{rk(E)}\chi(X,E\otimes L^k)\frac{\omega ^n}{n!} I_E}
   {\frac{\omega^n}{n!}}.
\nonumber
\eea
Recall that the moment map $\mu$ is given by
(\ref{moment}). Thus
\begin{eqnarray}
\dot{\mu}= \frac{d}{dt}[\sum_m(\frac{i}{2\pi}F+k\omega I)^m Td_{n-m}(X)]
= [\sum_m(\frac{i}{2\pi}F+k\omega I)^{m-1}\frac{i}{2\pi}\dot{F} Td_{n-m}(X)]
_{(sym)}.
\end{eqnarray}
Since $\dot{F}=\bar{\partial}\partial_H(\dot{h}h^{-1})$ for any
flow, we can combine these formulas and get
\bea
\dot{\mu}&=& [\sum_m(\frac{i}{2\pi}F+k\omega I)^{m-1}\frac{i}{2\pi}\bar{\partial}
\partial_H(\dot{h}h^{-1}) Td_{n-m}(X)]_{(sym)}\nonumber\\
           &=& [\sum_m(\frac{i}{2\pi}F+k\omega I)^{m-1}\frac{i}{2\pi}\times\nonumber\\
           &&\bar{\partial}\partial_H(\frac{-[\exp({\frac{i}{2 \pi}F + k \omega I_E})Td(X)]^{(2n)}
+ \frac{1}{rk(E)}\chi(X,E\otimes L^k)\frac{\omega ^n}{n!} I_E}
   {\frac{\omega^n}{n!}}) Td_{n-m}(X)]_{(sym)}\nonumber\\
           &=& [\sum_m(\frac{i}{2\pi}F+k\omega I)^{m-1}\frac{i}{2\pi}\bar{\partial}
\partial_H(\frac{\mu}{\frac{\omega ^n}{n!}}) Td_{n-m}(X)]_{(sym)}\nonumber\\
           &=& \sum_m \sum_{l=0}^{m-1}[\frac{i}{2\pi}C_l^{m-1}\omega ^l (\frac{i}{2\pi}F)^{m-1-l}
Td_{n-m}(X) \bar{\partial}\partial_H(\frac{\mu}{\frac{\omega ^n}{n!}})]_{(sym)}k^l.
\eea
This is the evolution equation for the moment map $\mu$. We observe that it is a
polynomial in $k$. The highest degree term has coefficient
$\Delta
\frac{\mu}{\frac{\omega ^n}{n!}}$ ,
but the rest of the terms involve F, which cannot be converted to $\mu$.
Q.E.D.

\section{Some explicit formulas}

The almost Hermitian-Einstein flow appears to be considerably more complicated than
the Hermitian-Einstein, or Donaldson heat flow. For the Donaldson heat
flow, as shown by Donaldson \cite{D},
important geometric quantities such as the curvature density
can be controlled, and the flow exists for all time.
But for the almost Hermitian-Einstein equation,
the flow for these
quantities are more complicated, and a full analysis is still unavailable.
In this section, we derive some explicit formulas to illustrate these
features. For the sake of simplicity, we consider only the case
of $X$ of dimension $2$ (when $X$ is dimension of 1, there is no difference between the almost Hermitian-Einstein flow and the Hermitian-Einstein flow) and $E$ of rank $1$, where the equation becomes already
complicated. In this case, the almost Hermitian-Einstein equation can be written as
\be
\frac{i}{2\pi}F\wedge \omega = \mu \frac{\omega ^2}{2}+\frac{1}{k}(\chi^2\frac{\omega^2}{2} - Td_X^2
- \frac{i}{2\pi}F\wedge Td_X^1 - (\frac{i}{2\pi}F)^2).
\ee
If we write down the local coordinate expression of the evolution equation
for the metric $h$ on the line bundle $E$, we get:
\be
\label{AHE21}
\dot{h}h^{-1} = - \hat{F} + \mu + \frac{1}{k} (\chi^2
+ \frac{1}{4\pi^2}(F_{i\bar{k}}F_{k\bar{i}} - \hat{F}^2)
+\frac{1}{8\pi^2}(F_{i\bar{k}}R_{k\bar{i}}- \hat{F}\hat{R})\nonumber\\
-\tilde{Td^2_X})\nonumber
\ee
where $F$ is the curvature of the bundle $E$, $R$ is the curvature of
the base manifold $X$, $\hat{F}$ means $F_{i\bar{i}}$, $\hat{R}$
is the scalar curvature and $\tilde{Td^2_X}$ is the coefficient of $\frac{\omega^2}{2}$ in the
second Todd class $Td^2$. It is known that $\tilde{Td^2_X}$ is a degree two polynomial in
curvature $R$. 
In the special case when the base metric $\omega$ is K\"ahler-Einstein, the flow reduces to 
\bea
\dot{h}h^{-1} &=& - \hat{F} + \mu + \frac{1}{k} (\chi^2
+ \frac{1}{4\pi^2}(F_{i\bar{k}}F_{k\bar{i}} - \hat{F}^2)-\tilde{Td^2_X}).\nonumber
\eea
Notice that the right hand side of the almost Hermitian-Einstein equation is non-linear in $F$. However, if we let $k\to\infty$, (\ref{AHE21}) goes to $$ \dot h h^{-1}= -\hat F +\mu,$$ the Hermitian-Einstein flow. 

Now, let's derive the evolution equation for $F$:
\begin{eqnarray}
\dot{F} &=& \bar\partial\partial_H(\dot{h}h^{-1}) \nonumber\\
&=& - \bar\partial\partial_H(\hat{F}-\mu-
    \frac{1}{k}(\chi^2 + \frac{1}{4\pi^2}(\hat{F}^2-F_{i\bar{k}}F_{k\bar{i}})
                + \frac{1}{8\pi^2}(\hat{F}\hat{R} - F_{i\bar{k}}R_{k\bar{i}})
                 +\tilde{Td^2_X}))  \nonumber\\
&=& - \bar\partial\partial_H(\hat{F}-
    \frac{1}{k}(\frac{1}{4\pi^2}(\hat{F}^2-F_{i\bar{k}}F_{k\bar{i}})
                + \frac{1}{8\pi^2}(\hat{F}\hat{R} - F_{i\bar{k}}R_{k\bar{i}})
                + \tilde{Td^2_X})) . \nonumber
\end{eqnarray}
Hence we have
\begin{eqnarray} 
\frac{d}{dt}|F|^2 &=& 2 Re(\dot{F},F) \nonumber\\
&=& 2 Re(- \bar\partial\partial_H(\hat{F}-
    \frac{1}{k}(\frac{1}{4\pi^2}(\hat{F}^2-F_{i\bar{k}}F_{k\bar{i}})
                + \frac{1}{8\pi^2}(\hat{F}\hat{R} - F_{i\bar{k}}R_{k\bar{i}})\nonumber\\
                &&+\tilde{Td^2_X})),F). \nonumber
\end{eqnarray}
Next we compute the $\frac{1}{k}$ - terms:
\begin{eqnarray}
(\bar{\partial}\partial_H(F_{i\bar{k}}R_{k\bar{i}}
)_{\bar{l}m}F_{l \bar{m}} &=& \nabla_{\bar{l}}\nabla_m(
F_{i\bar{k}}R_{k\bar{i}})F_{l \bar m} \nonumber\\
&=& (\nabla_{\bar{l}}\nabla_m
F_{i\bar{k}})R_{k\bar{i}}F_{l \bar m} + \nabla_m
F_{i\bar{k}} \nabla_{\bar l}R_{k\bar{i}}F_{l \bar m}+\nabla_{\bar l}
F_{i\bar{k}} \nabla_m R_{k\bar{i}}F_{l \bar m} \nonumber\\
&&+ F_{i\bar{k}} \nabla_{\bar l} \nabla_m R_{k \bar i} F_{l\bar m} \nonumber\\
&=& (\nabla_{\bar{l}}\nabla_i
F_{m\bar{k}})R_{k\bar{i}}F_{l \bar m} +  \nabla_m
F_{i\bar{k}} \nabla_{\bar l}R_{k\bar{i}}F_{l \bar m}+\nabla_{\bar l}
F_{i\bar{k}} \nabla_m R_{k\bar{i}}F_{l \bar m} \nonumber\\
&&+ F_{i\bar{k}} \nabla_{\bar l} \nabla_m R_{k \bar i} F_{l\bar m} \nonumber\\
&=& (\nabla_i\nabla_{\bar{l}}
F_{m\bar{k}}-F_{a\bar{k}}R^a_{mi\bar{l}}-F_{m\bar{a}}R^{\bar{a}}_{\bar k i \bar l })R_{k\bar{i}}F_{l \bar m}\nonumber \\
&&+ \nabla_m F_{i\bar{k}} \nabla_{\bar l}R_{k\bar{i}}F_{l \bar{m}}+\nabla_{\bar l}
F_{i\bar{k}} \nabla_m R_{k\bar{i}}F_{l\bar{m} } + F_{i\bar{k}} \nabla_{\bar l} \nabla_m R_{k \bar i} F_{l\bar m}\nonumber\\
&=& (\nabla_i\nabla_{\bar{k}}
F_{m\bar{l}})F_{l \bar{m}}R_{k\bar{i}} -F_{a\bar{k}}R^a_{mi\bar{l}}R_{k\bar{i}}F_{l \bar{m}}-F_{m\bar{a}}R^{\bar{a}}_{\bar k i \bar l  }
R_{k\bar{i}}F_{l \bar{m}}\nonumber \\
&&+ \nabla_m F_{i\bar{k}} \nabla_{\bar l}R_{k\bar{i}}F_{l \bar{m}}+\nabla_{\bar l}
F_{i\bar{k}} \nabla_m R_{k\bar{i}}F_{l \bar{m}}+ F_{i\bar{k}} \nabla_{\bar l} \nabla_m R_{k \bar i} F_{l\bar m}\nonumber\\
&=& \frac{1}{2}\nabla_i\nabla_{\bar{k}}|F_{m\bar{l}}|^2 R_{k\bar{i}}  -
\nabla_i F_{m\bar{l}}\nabla_{\bar{k}}F_{l\bar{m}}R_{k\bar{i}} \nonumber \\
&&+ \nabla_m F_{i\bar{k}} \nabla_{\bar l}R_{k\bar{i}}F_{l \bar{m}}+\nabla_{\bar l}
F_{i\bar{k}} \nabla_m R_{k\bar{i}}F_{l \bar{m}}\nonumber\\
&&-F_{a\bar{k}}R^a_{mi\bar{l}}R_{k\bar{i}}F_{l \bar{m}}-F_{m\bar{a}}R^{\bar{a}}_{\bar k i \bar l }
R_{k\bar{i}}F_{l \bar{m}}+ F_{i\bar{k}} \nabla_{\bar l} \nabla_m R_{k \bar i} F_{l\bar m}.\nonumber
\end{eqnarray}
Similarly
\begin{eqnarray}
(\bar{\partial}\partial_H(F_{i\bar{k}}F_{k\bar{i}}
)_{\bar{l}m}F_{l\bar{m}}
&=& \nabla_{\bar l}\nabla_{m}|F_{i\bar{k}}|^2 F_{l\bar{m}},\nonumber\\
(\bar{\partial}\partial_H(\hat{F}^2 )_{\bar{l}m}F_{l\bar{m}}
&=&
\hat{F}\triangle|F|^2 - 2|\nabla_i F_{m\bar l}|^2 \hat F+ 2 \nabla_m F_{i\bar{i}} \nabla_{\bar l}F_{j\bar{j}}F_{l \bar{m}}\nonumber\\
&& -2F_{a\bar{i}}R^a_{mi\bar{l}} \hat F F_{l \bar{m}}-2F_{m\bar{a}}R^{\bar{a}}_{\bar i i \bar l }
\hat F F_{l \bar{m}}, \nonumber\\
(\bar{\partial}\partial_H(\hat{F}\hat{R} )_{\bar{l}m}F_{l\bar{m}}
&=& \frac{1}{2}\hat{R}\triangle|F|^2 - |\nabla_i F_{m\bar l}|^2 \hat R + \nabla_m \hat F \nabla_{\bar{l}} \hat R F_{l\bar{m}}\nonumber \\
&&+ \nabla_{\bar l} \hat F \nabla_m \hat R F_{l \bar{m}}+\hat F \nabla_{\bar l} \nabla_m \hat RF_{l \bar{m}}\nonumber\\
&& -F_{a\bar{i}}R^a_{mi\bar{l}} \hat R F_{l \bar{m}}- F_{m\bar{a}}R^{\bar{a}}_{\bar i i \bar l }
\hat R F_{l \bar{m}}.\nonumber
\end{eqnarray}
If we follow Donaldson's argument, we can show when $rk(E)=1$, $\frac{d}{dt}|F|^2= -2Re( \bar\partial\partial_H\hat{F},F)$ implies an inequality 
\be
\label{HEine}
(\frac{d}{dt}-\Delta)|F|^2 \leq C|F|^2- 2|\nabla_i F_{m\bar l}|^2.
\ee
In our case the inequality becomes:
\begin{eqnarray}
    (\frac{d}{dt}-\triangle)|F|^2 &\leq&C|F|^2- 2|\nabla_i F_{m\bar l}|^2   -\frac{1}{k}\frac{1}{4\pi^2}(\hat{F}\triangle|F|^2-F_{k\bar{i}}\nabla_{i}\nabla_{\bar{k}}|F|^2)\nonumber\\   &-&\frac{1}{k}\frac{1}{16\pi^2}(\hat{R}\triangle|F|^2-R_{k\bar{i}}\nabla_{i}\nabla_{\bar{k}}|F|^2)  +\frac{1}{k}\frac{1}{8\pi^2} (
-\nabla_i F_{m\bar{l}}\nabla_{\bar{k}}F_{l\bar{m}}R_{k\bar{i}} \nonumber \\
&+& \nabla_m F_{i\bar{k}} \nabla_{\bar l}R_{k\bar{i}}F_{l \bar{m}}+\nabla_{\bar l}
F_{i\bar{k}} \nabla_m R_{k\bar{i}}F_{l \bar{m}}+C|F|^2 \nonumber \\
&+& 4|\nabla_i F_{m\bar l}|^2 \hat F- 4 \nabla_m \hat F \nabla_{\bar l}\hat F F_{l \bar{m}}+C|F|^3\nonumber\\
&+& |\nabla_i F_{m\bar l}|^2 \hat R - \nabla_m \hat F \nabla_{\bar{l}} \hat R F_{l\bar{m}}- \nabla_{\bar l} \hat F \nabla_m \hat R F_{l \bar{m}}+
C|F|^2+C|F|)\nonumber
\end{eqnarray}
We summarize this as,\\
{\bf Proposition}
{\it Let $X$ be of dimension 2 and $E$ be of rank 1. We have the following inequality
\begin{eqnarray}
&&\frac{d}{dt}|F|^2 + (\frac{1}{k}
\frac{1}{4\pi^2}\hat{F}+\frac{1}{k}\frac{1}{16\pi^2}\hat{R}-1)\triangle|F|^2
+\frac{1}{k}(-\frac{1}{4\pi^2}F_{k\bar{i}}-\frac{1}{16\pi^2}R_{k\bar{i}})\nabla_i\nabla_{\bar{k}}|F|^2)\nonumber\\
&\leq& C|F|^2- 2|\nabla_i F_{m\bar l}|^2 +\frac{1}{k}  ((\frac{1}{2\pi^2} \hat F+\frac{1}{8\pi^2}\hat R)|\nabla_i F_{m\bar l}|^2 - \frac{1}{2\pi^2} \nabla_m \hat F \nabla_{\bar l}\hat F F_{l \bar{m}}\nonumber\\
&-&\frac{1}{8\pi^2}\nabla_i F_{m\bar{l}}\nabla_{\bar{k}}F_{l\bar{m}}R_{k\bar{i}}  + \frac{1}{8\pi^2}\nabla_m F_{i\bar{k}} \nabla_{\bar l}R_{k\bar{i}}F_{l \bar{m}}+\frac{1}{8\pi^2}\nabla_{\bar l}
F_{i\bar{k}} \nabla_m R_{k\bar{i}}F_{l \bar{m}}\nonumber \\
&-&   \frac{1}{8\pi^2}\nabla_m \hat F \nabla_{\bar{l}} \hat R F_{l\bar{m}}- \frac{1}{8\pi^2}\nabla_{\bar l} \hat F \nabla_m \hat R F_{l \bar{m}}+C|F|^3+C|F|) .\nonumber
\end{eqnarray}
}
The inequality we get is much more complicated than the one, (\ref{HEine}), in the case of  Hermitian-Einstein flow. However, if we let $k\to\infty$, the inequality reduces to {(\ref{HEine})}
$$(\frac{d}{dt}-\Delta)|F|^2 \leq C|F|^2- 2|\nabla_i F_{m\bar l}|^2$$ as desired.
\\
If we assume that the base metric $g_{k \bar i}$ has a constant scalar curvature, the inequality simplifies to 
\begin{eqnarray}
&&\frac{d}{dt}|F|^2 + (\frac{1}{k}
\frac{1}{4\pi^2}\hat{F}+\frac{1}{k}\frac{1}{16\pi^2}\hat{R}-1)\triangle|F|^2
+\frac{1}{k}(-\frac{1}{4\pi^2}F_{k\bar{i}}-\frac{1}{16\pi^2}R_{k\bar{i}})\nabla_i\nabla_{\bar{k}}|F|^2)\nonumber\\
&\leq& C|F|^2- 2|\nabla_i F_{m\bar l}|^2 +\frac{1}{k}  ((\frac{1}{2\pi^2} \hat F+\frac{1}{8\pi^2}\hat R)|\nabla_i F_{m\bar l}|^2 - \frac{1}{2\pi^2} \nabla_m \hat F \nabla_{\bar l}\hat F F_{l \bar{m}}\nonumber\\
&-&\frac{1}{8\pi^2}\nabla_i F_{m\bar{l}}\nabla_{\bar{k}}F_{l\bar{m}}R_{k\bar{i}}  + \frac{1}{8\pi^2}\nabla_m F_{i\bar{k}} \nabla_{\bar l}R_{k\bar{i}}F_{l \bar{m}}+\frac{1}{8\pi^2}\nabla_{\bar l}
F_{i\bar{k}} \nabla_m R_{k\bar{i}}F_{l \bar{m}}\nonumber \\
&+&   C|F|^3+C|F|) .\nonumber
\end{eqnarray}
If we further assume that the base metric $g_{k \bar i}$ is K\"ahler-Einstein, the inequality simplifies to
\begin{eqnarray}
&&\frac{d}{dt}|F|^2 + (\frac{1}{k}
\frac{1}{4\pi^2}\hat{F}-1)\triangle|F|^2
+\frac{1}{k}(-\frac{1}{4\pi^2}F_{k\bar{i}})\nabla_i\nabla_{\bar{k}}|F|^2)\nonumber\\
&\leq& C|F|^2- 2|\nabla_i F_{m\bar l}|^2 +\frac{1}{k}  (\frac{1}{2\pi^2} |\nabla_i F_{m\bar l}|^2\hat F - \frac{1}{2\pi^2} \nabla_m \hat F \nabla_{\bar l}\hat F F_{l \bar{m}}\nonumber\\
&+&  C|F|^3+C|F|) .\nonumber
\end{eqnarray}

\end{document}